\documentclass[12pt,reqno]{amsart}

\usepackage{amssymb}
\usepackage[hypertex]{hyperref}

\begin{document}

\let\kappa=\varkappa
\let\eps=\varepsilon
\let\phi=\varphi

\def\Z{\mathbb Z}
\def\R{\mathbb R}
\def\C{\mathbb C}
\def\Q{\mathbb Q}

\def\OO{\mathcal O}
\def\CP{\C{\mathrm P}}
\def\RP{\R{\mathrm P}}

\def\conj{\overline}

\def\sest{\mathop{/\!\!/}}

\renewcommand{\Im}{{\mathop{\mathrm{Im}}\nolimits}}
\renewcommand{\Re}{{\mathop{\mathrm{Re}}\nolimits}}
\newcommand{\codim}{{\mathop{\mathrm{codim}}\nolimits}}

\newcommand{\demo}[1]{\noindent{\em #1}.}

\newtheorem{thm}{Theorem}
\newtheorem{lem}{Lemma}[section]
\newtheorem{prop}[lem]{Proposition}
\newtheorem{cor}[lem]{Corollary}

\theoremstyle{definition}
\newtheorem{df}[lem]{Definition}
\newtheorem{rem}[lem]{Remark}

\title{Levi problem and semistable quotients}
\author{Stefan Nemirovski}
\thanks{Supported in part by grants from DFG, RFBR, and the programme
``Leading scientific schools of Russia''}
\address{%
Steklov Mathematical Institute\hfill\break
\strut\hspace{8 true pt}
Mathematische Fakult\"at der Ruhr-Universit\"at Bochum}
\email{stefan@mi.ras.ru}

\maketitle

\bigskip

\section*{Introduction}
Given a complex space~$X$, the Levi problem asks for a global
characterisation of pseudoconvex unramified domains over~$X$.
For instance, if $p:U\to X$ is a pseudoconvex unramified domain
over a Stein manifold $X$, then $U$ is Stein by the classical
result of Docquier and Grauert~\cite{DG}.

The purpose of the present note is to introduce a broader class
of complex spaces within which the Levi problem will still
be tractable. The idea comes from the solution of the Levi problem
over Grassmanians given by Ueda~\cite{Ue}. In that work, the
realisation of the Grassmann manifold $\text{Gr}_\C(p,n)$
as the quotient of the space of $p\times n$ matrices of maximal rank
by the natural $\text{GL}(p,\C)$-action was used to show that
a pseudoconvex domain over $\text{Gr}_\C(p,n)$ is either Stein
or the whole of~$\text{Gr}_\C(p,n)$.

Let $G$ be a reductive complex Lie group.
A (normal) complex space $X$ is in class~${\mathcal Q}_G$
if it is biholomorphic to a semistable quotient (see \S\ref{QuotDef})
of the form
$$
X={\mathcal X}-\Sigma\sest G,
$$
where $\mathcal X$ is a Stein manifold with a holomorphic action
of $G$ and $\Sigma\subset\mathcal X$ is a $G$-invariant analytic subset
of codimension~$\ge 2$.

\begin{thm}
\label{main1}
If $p:U\to X$ is a pseudoconvex unramified domain
over $X\in{\mathcal Q}_G$, then $U\in {\mathcal Q}_G$.
\end{thm}

It is easy to show that every $X\in {\mathcal Q}_G$ is {\it holomorphically
precomplete}, that is, admits a holomorphic map $\iota:X\to X^\circ$
to a Stein complex space such that the induced homomorphism
$\iota^*:{\mathcal O}(X^\circ)\to{\mathcal O}(X)$ is an isomorphism
(see Proposition~\ref{Precomp}). If $X\in {\mathcal Q}_G$ is non-singular,
then the envelope of holomorphy of any domain over $X$ is pseudoconvex
and it follows that the domain itself is holomorphically precomplete.
The latter assertion can be proved directly for all $X\in {\mathcal Q}_G$.

\begin{thm}
\label{main2}
If $p:U\to X$ is an unramified domain over $X\in{\mathcal Q}_G$,
then $U$ is holomorphically precomplete.
\end{thm}

On the other hand, Grauert's examples~\cite{Gr} show that in general
a domain in a complex manifold need not be precomplete even if the
ambient manifold is projective algebraic. One can also exhibit
domains in Stein {\it spaces\/} that are not holomorphically
precomplete, see \cite{Bi} and~\cite{CD}.

Let us mention a couple of examples of complex spaces that belong
to (at least) one of the $\mathcal Q_G$'s. We have already seen
that the Grassmanian $\text{Gr}_\C(p,n)$ is in $\mathcal Q_{{\text{GL}(p,\C)}}$.
Somewhat more generally, the manifold of flags of type $0<p_1<\ldots<p_k<n$ in~$\C^n$
is in $\mathcal Q_{\text{GL}(p_1,\C)\times\text{GL}(p_2,\C)\times\cdots\times\text{GL}(p_k,\C)}$.
The latter observation was used by Adachi~\cite{Ad} to extend Ueda's
solution of the Levi problem to flag manifolds.

Another class of examples is given by toric varieties. By definition, a
toric variety is a normal variety~$X$ containing a Zariski open
subset isomorphic to the algebraic torus ${(\C^*)}^n$ such that
the natural action of the torus on itself extends to an action
on~$X$. The construction described by Cox~\cite{Co} shows that every
toric variety is a semistable quotient of $\C^N$ minus a finite
collection of linear subspaces by an action of a reductive
subgroup of~${(\C^*)}^N$. This property of toric varieties
and Ueda's approach were recently used by Ivashkovich to solve
the Levi problem on Hirzebruch surfaces~\cite[\S 4.6]{Iv}.

A useful observation is that if $X_1=\mathcal X_1-\Sigma_1\sest G_1$
and $X_2=\mathcal X_2-\Sigma_2\sest G_2$, then
$$
X_1\times X_2 =
\mathcal X_1\times \mathcal X_2 - (\Sigma_1\times \mathcal X_2)\cup (\mathcal X_1\times\Sigma_2)
\sest  G_1\times G_2.
$$
In other words, if $X_1\in\mathcal Q_{G_1}$ and $X_2\in\mathcal Q_{G_2}$,
then $X_1\times X_2\in\mathcal Q_{G_1\times G_2}$.
For instance, one can take a product with a Stein manifold considered
as the quotient of itself by the trivial action (cf.~\cite[\S 1.4]{Ne}).

\smallskip
\noindent
{\bf Convention.}
All complex spaces are assumed to be normal (and hence reduced
and locally irreducible).

\smallskip
\noindent
{\bf Note.} The present text is a slightly expanded and updated write-up of
the author's talk ``Quotients and pseudoconvexity'' given at the workshop
``Complex Geometry'' held at the Institut Henri Poincar\'e
on January 19--23, 2004.

\section{Holomorphic functions and geometric invariant theory}

\subsection{Semistable quotients and Stein spaces}
\label{QuotDef}
Suppose that a complex Lie group $G$ acts holomorphically on a complex space ${\mathcal X}$.
A {\it semistable quotient\/}%
\footnote{This terminology was introduced in~\cite{HMP}. An alternative term is
{\it analytic Hilbert quotient}.}
of ${\mathcal X}$ by this action is
a complex space $X$ such that there exists a holomorphic map
$\pi:{\mathcal X}\to X$ with the following properties:

\begin{itemize}
\item[1)] $\pi$ is $G$-invariant;
\item[2)] $\pi$ is locally Stein
(that is, each point in $X$ has a neighbourhood $U$ whose pre-image
$\pi^{-1}(U)$ is a Stein open set in ${\mathcal X}$);
\item[3)] ${\mathcal O}_X=\pi_*{\mathcal O}_{\mathcal X}^G$
(that is, the sheaf of holomorphic functions of $X$ is the sheaf of $G$-invariants
of the direct image sheaf of the structure sheaf of~${\mathcal X}$).
\end{itemize}
If a semistable quotient exists (which is not always the case), it
is unique up to biholomorphism and will be denoted
by~${\mathcal X}\sest G$.

\begin{prop}[{\rm Snow~\cite{Sn}}]
\label{SteinDown}
Let $G$ be a reductive group acting holomorphically on a Stein space~$\mathcal X$.
Then the semistable quotient ${\mathcal X}\sest G$ exists and is a Stein
space.
\end{prop}

The following proposition shows that, conversely, if the
semistable quotient is Stein, then the space $\mathcal X$ must be
Stein.

\begin{prop}[{\rm
Heinzner--Migliorini--Polito~\cite{HMP}}]
\label{SteinUp}
Let $G$ be a reductive group acting holomorphically on a complex space~$\mathcal X$.
Suppose that the semistable quotient ${\mathcal X}\sest G$ exists. If
${\mathcal X}\sest G$ is Stein, then $\mathcal X$ is also Stein.
\end{prop}

\begin{rem} Propositions~\ref{SteinDown} and~\ref{SteinUp} generalise
the result of Matsushima and Morimoto~\cite{MM} that the total space of a principle
bundle with a reductive structure group is Stein
if and only if the base of this bundle is Stein.
\end{rem}

\subsection{Local slices}
\label{slices}
Suppose that $G$ is a complex Lie group acting holomorphically on
a complex space~$\mathcal Y$. Assume that the orbit $G\cdot x$
of a point $x\in Y$ is closed in~$\mathcal Y$.
Denote by $G_x\subseteq G$ the isotropy subgroup of~$x$. A
{\it local slice\/} for the $G$-action at~$x$ is a locally closed
complex subspace $C$ of $\mathcal Y$ containing~$x$ such that
\begin{enumerate}
\item[1)] $B:= G_x\cdot C$ is a $G_x$-invariant Stein subspace of $\mathcal Y$;
\item[2)] $C$ is open in $B$;
\item[3)] the natural map of the twisted product $G\times_{G_x} B\to\mathcal Y$
is biholomorphic onto an open Stein subset of~$\mathcal Y$ that is saturated
with respect to the equivalence relation defined by $\mathcal O(\mathcal Y)^G$.
\end{enumerate}
The $G_x$-invariant subspace $B$ is called a (global) slice
for the $G$-action at~$x$.

\begin{prop}[{\rm Snow~\cite{Sn}}] Let $G$ be a reductive
group acting holomorphically on a Stein space~$\mathcal X$. For
any point $x\in\mathcal X$ with $G\cdot x$ closed and any
neighbourhood of~$x$, there exists a local slice for the
$G$-action at~$x$ contained in that neighbourhood.
\end{prop}

A geometric implication is that the points of the semistable quotient
of a complex space by a reductive group action are in one-to-one correspondence
with the closed orbits of the action. The projection to the quotient maps
the adherence of each closed orbit to the corresponding point in the quotient.

\subsection{Precompleteness and holomorphic convexity for quotients}
Let $\mathcal X$ be a Stein manifold with an action of a reductive group~$G$.
Assume that $\Sigma\subset\mathcal X$ is a $G$-invariant analytic
subset of codimension $\ge 2$ such that the semistable quotient
$X=\mathcal X-\Sigma\sest G$ exists.

By the universal property of the semistable quotient, the inclusion
$$
\iota:\mathcal X - \Sigma \hookrightarrow \mathcal X
$$
induces a canonical holomorphic map
$$
\iota_G:\mathcal X-\Sigma\sest G\longrightarrow \mathcal X\sest G,
$$
where the latter space is Stein by Snow's result. Since $\codim_\C \Sigma\ge2$,
the map $\iota$ induces an isomorphism of the algebras of holomorphic
functions and hence so does~$\iota_G$. We have thus proved the following:

\begin{prop}
\label{Precomp}
Any quotient of the form $X=\mathcal X-\Sigma\sest G$ is holomorphically
precomplete.
\end{prop}

More can be said about $X$ if more is known
about the map $\iota_G$.

\begin{prop}
The quotient $X=\mathcal X-\Sigma\sest G$ is holomorphically convex
if and only if the map $\iota_G$ is proper.
\end{prop}

\demo{Proof} Recall that $X$ is holomorphically
convex if and only if there exists a Remmert reduction map, that
is, a {\it proper\/} map $\rho:X\to Y$ onto a Stein space such that
$\rho^*:{\mathcal O}(Y)\to{\mathcal O}(X)$ is an isomorphism. It
follows that any holomorphic map from $X$ to a Stein space factors
through~$\rho$.

Now, if $\iota_G$ is proper, then it is the required reduction map.
Conversely, if a reduction map $\rho$ exists, then $\iota_G$ factors through
it so that $\iota_G=\varphi\circ\rho$. The map $\phi$ is a map of Stein spaces
inducing an isomorphism on the algebras of holomorphic functions and hence
a biholomorphism. \qed

\subsection{Quotients without non-constant functions}
Suppose that $X=\mathcal X -\Sigma\sest G$ is such that $\mathcal O(X)=\C$.
(For instance, $X$ may be compact and connected.) Then
$\mathcal O(\mathcal X\sest G)=\C$ as well and hence
the Stein complex space $\mathcal X\sest G$ is a singleton.
As we have seen in \S\ref{slices}, this implies that
the $G$-action on $\mathcal X$ has a single closed orbit.
Since $\mathcal X$ is non-singular,
another result of Snow~\cite[Corollary 5.6(3)]{Sn}
shows that $\mathcal X$ is $G$-equivariantly biholomorphic
to a homogeneous vector bundle over the homogeneous affine
variety $G/H$ for a reductive complex subgroup $H\subseteq G$
and each fibre of this bundle is $H$-equivariantly biholomorphic
to a rational representation space of $H$.

Thus, it follows from the definitions that
$$
X = \C^n - \Xi\sest H,
$$
where $H$ acts on $\C^n$ by a rational linear representation with the origin $\{0\}\subset\C^n$
as the only closed orbit. The $H$-invariant analytic subset $\Xi\subset\C^n$
corresponds to the intersection of $\Sigma\subset\mathcal X$
with a fibre of the bundle~$\mathcal X\to G/H$.

\section{Pseudoconvex domains}
\subsection{Unramified domains}
An unramified domain over a complex space $X$ is a pair
$(U,p)$ consisting of a connected Hausdorff topological space $U$
and a locally homeomorphic map $p:U\to X$.
A domain is {\it schlicht\/} if the map $p$
is injective.
There exists a unique complex structure on the space~$U$
such that the projection $p:U\to X$ is a locally biholomorphic map.

\begin{df}
An unramified domain $p:U\to X$ is said to be {\it pseudoconvex\/}%
\footnote{In the literature, such domains are also called
{\it locally pseudoconvex\/} or {\it  locally Stein}.}
if for every point $x\in X$
there exists a neighbourhood $V\ni x$ such that its
pre-image $p^{-1}(V)$ is a Stein open subset of~$U$.
\end{df}

A domain $(U_1,p_1)$ is contained in $(U_2,p_2)$
if there exists a map $j:U_1\to U_2$ such that
$p_2\cdot j=p_1$. Notice that the map $j$ is {\it a~priori\/}
only locally biholomorphic. However, if $(U_1,p_1)$ is a
schlicht domain, then $j$ is a set-theoretic injection as well.

\subsection{Boundary points}
Let $p:U\to X$ be an unramified domain over a complex space~$X$.
Following Grauert and Remmert~\cite[Definition 4]{GR}, we introduce the notion
of a {\it boundary point of\/}~$(U,p)$. Namely, a boundary point
is a filtre~${\mathfrak V}$ of open connected subsets $V\subset U$
such that
\begin{itemize}
\item[1)] the filtre $\{p(V) \mid V\in\mathfrak V\}$
converges to a point~$x_0\in X$;
\item[2)] for every neighbourhood $W\ni x_0$,
the filtre $\mathfrak V$ contains precisely one connected
component of $p^{-1}(W)$;
\item[3)] the filtre $\mathfrak V$ has no accumulation point in~$U$.
\end{itemize}
The set of boundary points is called the boundary~$\partial U$ of
the domain $(U,p)$. The union $\widehat U=U\cup\partial U$ has a
structure of a Hausdorff topological space such that the map
$p:U\to X$ extends to a continuous map $\widehat p:\widehat U\to
X$. Note that the extension $\widehat p$ is not necessarily
locally homeomorphic.

A domain $(U,p)$ is called {\it pseudoconvex at a boundary
point\/} $a\in\partial U$ if there exists a neighbourhood $V\ni a$
in~$\widehat U$ such that $V\cap U$ is a Stein space. It follows
from Oka's theorem~\cite{O} that a domain over a complex manifold
is pseudoconvex if and only if it is pseudoconvex
at every boundary point.

\subsection{Unramified extension}
Let now $p:U\to X$ be a domain over a complex manifold~$X$.
Suppose that $S\subset X$ is an analytic subset of positive
codimension. A boundary point $a\in\partial U$ is called {\it
removable along\/}~$S$ if there exists a neighbourhood $V\ni a$
such that $(V,\widehat p)$ is a schlicht domain over~$X$ and the
set $V\cap\partial U$ is contained in $\widehat p^{-1}(S)$. If
$R\subset\partial U$ is the set of removable points, then
$(U^*,p^*)=(U\cup R, \widehat p|_{U\cup R})$ is an unramified
domain over $X$ that is called the {\it extension of $U$ along~$S$}.

The following result was proved by Ueda~\cite[p.\,393]{Ue} as a corollary
of~\cite[Satz~4]{GR}.

\begin{prop} \label{Ext}
Let $(U,p)$ be an unramified domain over a
complex manifold~$X$, and let $S$ be an analytic subset of
positive codimension in~$X$. Assume that $(U,p)$ is pseudoconvex
at every boundary point lying over~$X-S$. Then the extension of
$(U,p)$ along $S$ is pseudoconvex.
\end{prop}

In particular, if there exists no boundary point that is removable along~$S$,
then the domain $(U,p)$ is itself pseudoconvex.

\begin{rem}
Examples from Grauert's paper~\cite{Gr} show that Proposition~\ref{Ext}
can be false over a singular space $X$ even if the singularities
are isolated and the domain is schlicht.
A~simple non-schlicht example is given in \S\ref{exampleZ2} below.
\end{rem}

\subsection{Proof of Theorem~\ref{main1}}
\label{ProofMain1}
Assume now that $\mathcal X$ is a Stein manifold and consider
a pseudoconvex domain $(U,p)$ over a semistable quotient
$X={\mathcal X}-\Sigma\sest G$. Let $\pi:{\mathcal X}-\Sigma\to X$
be the natural projection.

There is a canonical `pull-back' of $(U,p)$ to a $G$-equivariant
unramified domain $(\mathcal U,\wp)$ over $\mathcal X$.
Namely, let
$$
\mathcal U=\{(u,x)\in U\times\mathcal X-\Sigma\mid p(u)=\pi(x)\}
$$
be the fibred product of $p$ and $\pi$. The space $\mathcal U$
fits into the commutative diagram
$$
\begin{array}{ccc}
\mathcal U&\stackrel{\wp}{\longrightarrow}&\mathcal X - \Sigma\\
{\varpi}\downarrow\phantom\pi&&\downarrow\lefteqn{\pi}\\
U&\stackrel{p}{\longrightarrow}& X
\end {array}
$$
where the maps $\wp$ and $\varpi$ are obvious
projections sending a pair $(u,x)\in\mathcal U$ to $x$ and~$u$,
respectively. The space $\mathcal U$ inherits a natural complex
structure such that $\varpi$ is a semistable quotient map onto~$U$
and the map $\wp$ is locally biholomorphic. It follows that
$\mathcal U$ is a disjoint union of (finitely many) unramified
domains over $\mathcal X$.

Proposition~\ref{SteinUp} shows that $(\mathcal U,\wp)$ is pseudoconvex
at every boundary point over $\mathcal X - \Sigma$. Hence, its unramified extension
$(\mathcal U^*, \wp^*)$ along $\Sigma$ is a finite union of pseudoconvex domains
over $\mathcal X$ by Proposition~\ref{Ext}.
By the Docquier--Grauert theorem~\cite{DG}, the manifold $\mathcal U^*$ is Stein.
Thus, we obtain the following slightly more precise version of Theorem~\ref{main1}.

\begin{prop}
\label{Quot}
Let $p:U\to X$ be a pseudoconvex domain over a semistable quotient
$X=\mathcal X-\Sigma\sest G$. Then $U=\mathcal U^* - (\wp^*)^{-1}(\Sigma)\sest G$, where
$\wp^*:\mathcal U^*\to\mathcal X$ is a $G$-equivariant Stein manifold
spread over $\mathcal X$.
\end{prop}

The manifold $\mathcal U^*$ may be disconnected if the group $G$ is disconnected
(see~\S\ref{exampleZ2}). In that case, however, $U=\mathcal U^*_0-(\wp^*)^{-1}(\Sigma)\sest G_0$,
where $\mathcal U^*_0$ is a connected component of $\mathcal U^*$ and $G_0$ is
the subgroup of $G$ leaving $\mathcal U^*_0$ invariant.

\begin{rem}
Suppose that the quotient $X=\mathcal X -\Sigma\sest G$ is Stein.
Then $\mathcal X-\Sigma$ is Stein by Proposition~\ref{SteinUp} and therefore $\Sigma=\varnothing$.
In this situation, the proof of Proposition~\ref{Quot} shows that any pseudoconvex
domain over $X$ is in fact Stein. The point here is that $X$ may have non-isolated
singularities and so the generalisation of the Docquier--Grauert theorem to Stein
spaces with isolated singularities due to Col\c{t}oiu and Diederich~\cite{CD2}
cannot be applied directly.
\end{rem}

\begin{rem}
The argument from the proof of Theorem~\ref{main1} can be used
to solve the Levi problem over $X=\mathcal X -\Sigma\sest G$ in rather
explicit terms provided that the action of $G$ on $\mathcal X$ is sufficiently
well understood. Examples can be found in the aforementioned
papers of Ueda~\cite{Ue,Ue2}, Adachi~\cite{Ad}, and Ivashkovich~\cite{Iv}.
\end{rem}

\section{Precompleteness of unramified domains}

\subsection{Envelopes of holomorphy}
For an unramified domain $p:U\to X$ over a complex space, its {\it envelope of holomorphy\/}
is a domain $\widetilde p:\widetilde U\to X$ with a locally biholomorphic
map $j:U\to\widetilde U$ such that
\begin{itemize}
\item[1)] $\widetilde p\circ j=p\;$;
\item[2)] for every holomorphic function $f\in\mathcal O(U)$ there exists
a holomorphic function $\widetilde f\in\mathcal O(\widetilde U)$ such that
$\widetilde f\circ j=f$;
\item[3)] if a domain $q:V\to X$ and a locally biholomorphic map $k:U\to V$
satisfy (1) and (2), then there exists a locally biholomorphic map
$\ell:V\to\widetilde U$ such that $\ell\circ k=j$ and $\widetilde p\circ \ell=q$.
\end{itemize}
In other words, $(\widetilde U, \widetilde p)$ is the maximal domain over $X$
containing $(U,p)$ such that all holomorphic functions from $U$ can be
holomorphically extended to $\widetilde U$.

The envelope of holomorphy exists for every domain and is unique up to a
natural isomorphism. If $X$ is non-singular, then the envelope is pseudoconvex
by the Cartan--Thullen--Oka theorem. It follows that if $X$ is a complex manifold
in $\mathcal Q_G$, then the envelope of holomorphy of any domain over~$X$ is
also in $\mathcal Q_G$ and hence the domain is holomorphically precomplete.
It turns out that this corollary (but {\it not\/} the Cartan--Thullen--Oka theorem!)
can be extended to singular complex spaces in~$\mathcal Q_G$.

\begin{lem}
Let $p:U\to \mathcal X$ and $q:V\to \mathcal X$ be two domains over a Stein manifold
and $j:U\to\widetilde U$ and $k:V\to\widetilde V$ the maps into their
respective envelopes of holomorphy. Then for every biholomorphic map
$f:U\to V$ there exists a unique biholomorphic extension
$\widetilde f:\widetilde U\to\widetilde V$
such that $\widetilde f\circ j(x)=k\circ f(x)$ for all $x\in U$.
\end{lem}

\demo{Proof} The lemma is well-known in the case $\mathcal X=\mathbb C^n$ (see,
for instance, \cite{Ma}). The following proof in the general case is essentially
contained in~\cite{Ke}.

The map $k\circ f:U\to\widetilde V$ extends to a holomorphic map $\widetilde f:\widetilde U\to \widetilde V$.
Indeed, since $\widetilde V$ is Stein  (because $\mathcal X$ is a Stein manifold and so the theorems
of Cartan--Thullen--Oka and Docquier--Grauert apply), there exists a proper holomorphic embedding
$\iota:\widetilde V\to\mathbb C^N$. By the definition of the enevelope, the map
$\iota\circ k\circ f:U\to\C^N$ extends to a holomorphic map $F:\widetilde U\to\C^N$.
The image of $F$ is contained in $\iota(\widetilde V)$ by the uniqueness theorem
and hence we can define $\widetilde f:=\iota^{-1}\circ F$.

Applying the same argument with $f^{-1}$ instead of $f$, we see that the map
$\widetilde f$ is biholomorphic.\qed

\begin{cor}
\label{EquivEnv}
Suppose that a group $G$ acts on a Stein manifold $\mathcal X$ by biholomorphisms.
Let $p:U\to \mathcal X$ be a $G$-equivariant union of unramified domains over $\mathcal X$.
Then the union of envelopes of holomorphy of these domains is
$G$-equivariant.
\end{cor}

\demo{Proof} Applying the preceding lemma to every element of $G$, we obtain
a $G$-action on the union of the envelopes. The compatibility of this
action with the $G$-action on $\mathcal X$ follows by the uniqueness theorem.\qed

\subsection{Proof of Theorem~\ref{main2}}
Suppose that $X=\mathcal X - \Sigma\sest G$ and let $p:U\to X$ be an unramified
domain over~$X$. Let $\wp:\mathcal U\to \mathcal X$ be its pull-back
to $\mathcal X$ (see the proof of Proposition~\ref{Quot} above) and consider its component-wise
envelope of holomorphy $(\widetilde{\mathcal U},\widetilde{\wp})$.
It follows from Corollary~\ref{EquivEnv} that
$\widetilde{\mathcal U}$ is a Stein manifold with a holomorphic $G$-action.
Furthermore, since $G$-invariance is preserved by analytic continuation,
the space of $G$-invariant functions on $\widetilde{\mathcal U}$ is isomorphic
to $\mathcal O(U)$.

The quotient $U^\circ=\widetilde{\mathcal U}\sest G$ exists and is a Stein space
by Proposition~\ref{SteinDown}. The map $U\to U^\circ$ induced by the natural map of
$\mathcal U$ into its envelope of holomorphy induces an isomorphism on
the algebras of holomorphic functions by the preceding remark. Hence,
$U$ is holomorphically precomplete as claimed.\qed

\begin{rem} In general, $U^\circ$ cannot be represented as a domain over~$X$.
Even if the projection $p:U\to X$ extends to a holomorphic map from $U^\circ$
to $X$, this map may be ramified. The latter situation may occur, for instance,
if $X$ is a {\it singular\/} Stein space.
\end{rem}

\subsection{Example}
\label{exampleZ2}
Consider the obvious action of
the group $G=\Z/2\Z$ on ${\mathcal X}=\C^2$ given by $(z_1,z_2)\mapsto(-z_1,-z_2)$.
The semistable quotient $X=\mathcal X\sest G$ is a Stein (in fact, affine) complex surface
with an $A_1$ singularity. Let $o\in X$ be the singular point. Then the restriction of
the quotient map $\pi$ to $\C^2-\{0\}=\pi^{-1}(X-\{o\})$ defines an unramified domain $(U,p)$
over~$X$. It is easy to see that this is a domain of holomorphy over $X$. However,
it is not pseudoconvex and does not admit a pseudoconvex unramified extension along~$\{o\}$.

Applying the construction from the proof of Theorem~\ref{main2} gives the following result. Firstly,
the pull-back of $(U,p)$ to ${\mathcal X}=\C^2$ is the union of two schlicht domains
$\mathop{\mathrm{id}}:\C^2-\{0\}\to\C^2-\{0\}$ with the $G$-action taking the point $(z_1,z_2)$ on one domain to the
point $(-z_1,-z_2)$ on the other. The component-wise envelope of holomorphy consists of two
schlicht domains $\mathop{\mathrm{id}}:\C^2\to\C^2$ with the extended `intertwining' $G$-action. The quotient
by this action is biholomorphic to~$\C^2$, which is the natural holomorphic completion of
$U=\C^2-\{0\}$. Notice that the projection $p:U\to X$ does indeed extend to $U^\circ=\C^2$ but this
extension is a ramified map.

\end{document}